\documentclass[a4paper, 10pt]{amsart}
\usepackage{filecontents}

\usepackage{etoolbox}
\makeatletter
\let\ams@starttoc\@starttoc
\makeatother
\usepackage[parfill]{parskip}
\makeatletter
\let\@starttoc\ams@starttoc
\patchcmd{\@starttoc}{\makeatletter}{\makeatletter\parskip\z@}{}{}
\makeatother

\usepackage[colorlinks=true,linkcolor=blue,citecolor=blue,urlcolor=blue]{hyperref}
\usepackage{amsthm}
\usepackage{bookmark}
\usepackage{thmtools}

\usepackage{amssymb}
\usepackage{amsmath}
\usepackage{fancyhdr}
\usepackage{esint}
\usepackage{enumerate}

\usepackage{pictexwd,dcpic}

\swapnumbers
\declaretheorem[name=Theorem,numberwithin=section]{thm}
\declaretheorem[name=Remark,style=remark,sibling=thm]{rem}
\declaretheorem[name=Lemma,sibling=thm]{lemma}
\declaretheorem[name=Proposition,sibling=thm]{prop}
\declaretheorem[name=Definition,style=definition,sibling=thm]{defn}
\declaretheorem[name=Corollary,sibling=thm]{cor}
\declaretheorem[name=Assumption,style=remark,sibling=thm]{ass}
\declaretheorem[name=Example,style=remark,sibling=thm]{example}

\numberwithin{equation}{section}

\usepackage{cleveref}
\crefname{lemma}{Lemma}{Lemmata}
\crefname{prop}{Proposition}{Propositions}
\crefname{thm}{Theorem}{Theorems}
\crefname{cor}{Corollary}{Corollaries}
\crefname{defn}{Definition}{Definitions}
\crefname{example}{Example}{Examples}
\crefname{rem}{Remark}{Remarks}
\crefname{ass}{Assumption}{Assumptions}

\renewcommand{\(}{\left(}
\renewcommand{\)}{\right)}

\renewcommand{\-}{\bar}

\newcommand{\cn}{\colon}

\newcommand{\R}{\mathbb{R}}

\renewcommand{\S}{\mathbb{S}}

\newcommand{\Di}{\mathbb{D}}

\renewcommand{\a}{\alpha}
\renewcommand{\b}{\beta}
\newcommand{\g}{\gamma}
\renewcommand{\d}{\delta}
\newcommand{\e}{\epsilon}
\renewcommand{\k}{\kappa}
\renewcommand{\l}{\lambda}

\renewcommand{\o}{\omega}
\newcommand{\D}{\Delta}

\newcommand{\s}{\sigma}
\newcommand{\p}{\varphi}

\newcommand{\z}{\zeta}

\newcommand{\del}{\partial}

\newcommand{\inpr}[2]{\langle #1,#2 \rangle}
\newcommand{\fr}[2]{\frac{#1}{#2}}

\newcommand{\Thm}{\begin{thm}}
\newcommand{\eThm}{\end{thm}}
\newcommand{\Def}{\begin{defn}}
\newcommand{\eDef}{\end{defn}}
\newcommand{\Prop}{\begin{prop}}
\newcommand{\eProp}{\end{prop}}
\newcommand{\Rem}{\begin{rem}}
\newcommand{\eRem}{\end{rem}}
\newcommand{\Lem}{\begin{lemma}}
\newcommand{\eLem}{\end{lemma}}
\newcommand{\eq}{\begin{equation}}
\newcommand{\eeq}{\end{equation}}
\newcommand{\Ex}{\begin{example}}
\newcommand{\eEx}{\end{example}}
\newcommand{\pf}{\begin{proof}}
\newcommand{\epf}{\end{proof}}
\newcommand{\Cor}{\begin{cor}}
\newcommand{\eCor}{\end{cor}}
\newcommand{\Ass}{\begin{ass}}
\newcommand{\eAss}{\end{ass}}
\newcommand{\SAl}{\begin{align}\begin{split}}

\newcommand{\PBVP}[9]{\begin{align}\begin{split}\label{PBVP} #1=#2\ &\mathrm{in}\ #3\\
										#4=#5\ &\mathrm{on}\ #6\\
										#7=#8\ &\mathrm{on}\ #9\end{split}\end{align}}

\newcommand{\ra}{\rightarrow}
\newcommand{\hra}{\hookrightarrow}
\newcommand{\mt}{\mapsto}

\newcommand{\mc}{\mathcal}

\newcommand{\mrm}{\mathrm}

\newcommand{\hp}{\hphantom}

\protected\def\ignorethis#1\endignorethis{}
\let\endignorethis\relax

\begin{document}

\title{The inverse mean curvature flow perpendicular to the sphere}
\author{Ben Lambert and Julian Scheuer}

\begin{abstract}
We consider the smooth inverse mean curvature flow of strictly convex hypersurfaces with boundary embedded in $\R^{n+1},$ which are perpendicular to the unit sphere from the inside. We prove that the flow hypersurfaces converge to the embedding of a flat disk in the norm of $C^{1,\b},$ $\b<1.$
\end{abstract}
\date{\today}
\address{Ben Lambert, University of Konstanz, Zukunftskolleg, Box 216, 78457 Konstanz, Germany}
\email{benjamin.lambert@uni-konstanz.de}
\address{Julian Scheuer, Ruprecht-Karls-Universit{\"a}t, Institut f{\"u}r Angewandte Mathematik, Im Neuenheimer Feld 294,
69120 Heidelberg, Germany}
\email{scheuer@math.uni-heidelberg.de}
\thanks{The first author is a ZIF-Marie Curie fellow at the Zukunftskolleg, Konstanz.}
\thanks{The second author is being supported by the DFG}

\maketitle
\tableofcontents

\section{Introduction}
We consider the inverse mean curvature flow in $\R^{n+1}$ with a Neumann boundary condition in a sphere. Let $\Di=\Di^{n}$ be the $n$-dimensional unit disk and $\~N$ be the outward unit normal of the inclusion $\S^{n}\hra\R^{n+1}.$ Then we consider a family of embeddings 
\eq X\cn [0,T^{*})\times\Di\hra \R^{n+1}\eeq with a normal vector field $N,$ the choice of which will be specified in a natural manner later, 
such that
\begin{subequations}\label{Flow}\begin{align}\dot{X}&=\fr{1}{H}N, \label{Floweq}\\
		X(\del \Di)&=\del X(\Di)\subset\label{Flow2} \S^{n}, \\
		0&=\inpr{N_{|\del\Di}}{\~N(X_{|\del\Di})}, \label{BoundCond}\\
		\inpr{\dot{\g}(0)}{\~N}&\geq 0\quad\forall \g\in C^{1}((-\e,0],M_{t})\cn \g(0)\in\del X(\Di). \label{Flow4}\end{align}\end{subequations}
We prove the following result.

\Thm\label{Main}
Let \eq X_{0}\cn \Di\hra M_{0}\subset\R^{n+1}\eeq be the embedding of a smooth and strictly convex hypersurface with normal vector field $N_{0},$ such that 
\begin{subequations}\label{initial}\begin{align} \inpr{\dot{\g}(0)}{\~N}&\geq 0\quad\forall \g\in C^{1}((-\e,0],M_{0})\cn \g(0)\in\del X_{0}(\Di),\label{initial1}\\
									 X_{0}(\del \Di)&\subset\S^{n},\\
									 \inpr{N_{0|\del\Di}}{\~N_{|\del\Di}}&=0.\label{initial3}\end{align}\end{subequations}
Then there exists a finite time $T^{*}<\infty,$ $\alpha>0$ and a unique solution
\eq X\in C^{1+\fr{\a}{2},2+\a}([0,T^{*})\times\Di)\cap C^{\infty}((0,T^{*})\times \Di,\R^{n+1})\eeq of \eqref{Flow} with initial hypersurface $M_{0},$ such that the embeddings $X_{t}$ converge to the embedding of a flat unit disk as $t\ra T^{*},$ in the sense that the height of the $M_{t}=X(t,M)$ over this disk converges to $0.$ 
%in the norm of $C^{1,\b},$ $0<\b<1.$ 
\eThm

\Rem
The norm of convergence of the $M_{t}$ to the disk will be specified in \cref{Convergence}, when we will have developed a suitable coordinate system to describe the $M_{t}.$
\eRem

Our motivation for treating this problem arises from several directions. First of all, the inverse mean curvature flow (IMCF) has proven to be a useful tool in the theory of geometric inequalities, cf.~\cite{HuiskenIlmanen:/2001} for the probably most famous result in this direction. The works which describe the asymptotic behavior of the IMCF in Euclidean space include \cite{Gerhardt:/1990} and \cite{Urbas:/1990}, whereas in the hyperbolic space we refer to \cite{Gerhardt:11/2011} and \cite{Scheuer:06/2014}. Those works deal with closed hypersurfaces. 

Few years ago, the Ph.D. thesis \cite{Marquardt:/2012} written by Thomas Marquardt appeared, also cf.~\cite{Marquardt:07/2013}. Here the IMCF of hypersurfaces with boundary was considered and the embedded flowing hypersurfaces were supposed to be perpendicular to a convex cone in $\R^{n+1}.$ However, short-time existence was derived in a much more general situation, in other ambient spaces and with other supporting hypersurfaces besides the cone. It appears to be a natural question, whether one can also obtain nice convergence results if one imposes perpendicularity to other hypersurfaces. Inspired by a recent result about rigidity of hypersurfaces in the sphere by Matthias Makowski and the second author, cf.~\cite{MakowskiScheuer:/2013}, Oliver Schn{\"u}rer suggested to the authors that this rigidity result might be helpful to consider the IMCF for hypersurfaces which are perpendicular to the sphere. Indeed, we were able to prove his conjecture that this flow must drive strictly convex hypersurfaces into the embedding of a disk.

The equivalent problem for the mean curvature flow was treated by Axel Stahl in \cite{Stahl:06/1996} and \cite{Stahl:08/1996}, in which the flow was shown to contract to a point. Other choices of boundary manifolds for a graphical mean curvature flow have shown convergence of the flow to flat disks, see for example \cite{Lambert:/2014} and \cite{Huisken:02/1989}, as well as \cite{GigaOhnumaSato:05/1999} for a levelset approach.

The proof of \cref{Main} is ordered as follows: In \cref{Notation} we agree on notation and in \cref{Eveq} we collect the relevant evolution equations and boundary derivatives. In \cref{Heightest} we make height and gradient estimates for convex hypersurfaces perpendicular to the sphere, which is of interest independently. In particular there follows that if the boundary of a convex manifold is contained in a hemisphere, then we have a lower height bound on the manifold. In \cref{MoebiusCoord} we show that the flow may be written graphically. In \cref{CurvEst} we use the results of \cref{Heightest} to demonstrate the two key estimates which in conjunction with rigidity results of \cite{MakowskiScheuer:/2013} give the theorem. The first of these is that while the the boundary stays away from an equator, a convex flow has a lower bound on $H$. The second shows that the flow remains convex up until the singular time. Therefore, due to rigidity at the boundary, $\partial M$ must flow to an equator and so $M$ must flow to a flat disk assuming that the flow may be suitably extended. In \cref{Flatdisk} we clarify the necessary PDE existence results and show $C^{1, \beta}$ convergence.

\section{Setting and notation}\label{Notation}
There are various embeddings involved in \eqref{Flow}, namely the inclusion
\eq x\cn \S^{n}\hra \R^{n+1},\eeq
the flow embeddings of the form
\eq X\cn \Di\hra\R^{n+1},\eeq the inclusion
\eq z\cn \del\Di\hra \Di,\eeq
as well as the derived embedding
\eq\label{IndEmb1} y\cn \del\Di \hra \S^{n}\eeq
satisfying
\eq\label{IndEmb2} X\circ z=x\circ y.\eeq
Throughout this paper, we stick to the coordinate based notation for tensors.

 Geometric quantities in $\R^{n+1}$ are denoted by a bar, e.g. $(\-g_{\a\b})$ for the Euclidean metric, where greek indices range from $0$ to $n.$ We will also write $\inpr{\cdot}{\cdot}$ for the Euclidean scalar product.
 
 Geometric quantities in $\S^{n}$ are denoted by a check, e.g. $(\check g_{ij})$ for the induced metric of the embedding $x,$ where latin indices range from $1$ to $n.$
 
 Induced quantities of embeddings $\Di\hra\R^{n+1}$ are denoted by latin letters, e.g. the embeddings $X$ induce metrics $(g_{ij}),$ normal vector fields $N$ and a second fundamental forms $(h_{ij}),$ such that we have the Gaussian formula
 \eq X_{ij}^{\a}=-h_{ij}N^{\a}.\eeq
A hypersurface $M\hra\R^{n+1}$ is called strictly convex, if $N$ can smoothly be chosen, such that $(h_{ij})$ is positive definite. For a strictly convex hypersurface we will choose $N$ like this.
 
 Induced quantities of embeddings to $\del\Di\hra\S^{n}$ are denoted by greek letters, e.g. the embeddings $y$ induce metrics $(\g_{IJ}),$
 normal vector fields $\nu$ and second fundamental forms $(\eta_{IJ}),$ where capital latin indices range from $2$ to $n.$ 

Coordinate systems in $\del\Di$ will be denoted by $(\xi^{I}),$ $2\leq I\leq n.$

Define $H$ to be the mean curvature of the embeddings $X,$
\eq H=g^{ij}h_{ij},\eeq
where $(g^{ij})$ is the inverse of $(g_{ij}).$ 

For an embedded manifold $M^{n}\hra N^{n+1}$ and a function $u\cn M\ra\R,$ covariant derivatives with respect to the induced metric are denoted by indices, e.g. $u_{ij}.$
If ambiguities are possible, e.g. in the case of tensor derivation, covariant derivatives are denoted by a semicolon, e.g. $h_{ij;k}.$ Standard partial derivatives are denoted by a comma, e.g. $u_{i,j}.$

\section{Evolution equations and boundary derivatives}\label{Eveq}
For the inverse mean curvature flow the interior evolution equations are well-known. We need the spatial boundary derivatives of various curvature quantities, when the supporting hypersurface is a sphere. The calculations are quite similar to those in \cite{Marquardt:/2012} and \cite{Stahl:08/1996}. For the sake of completeness and for a better comprehensibility of the different notation, let us derive them in detail.

\Rem\label{Shorttime}
Short-time existence for the flow \eqref{Flow} was derived in \cite[Thm.~2.12]{Marquardt:/2012}.
Thus we are justified to use \eqref{Flow} to calculate the boundary derivatives.
\eRem

\Rem\label{gsplit}
Due to \eqref{BoundCond} we obtain that
\eq \~N\in X_{*}(T\Di)\eeq
and thus at boundary points there holds
\eq \~n\equiv(\inpr{X_{k}}{\~N})\in T^{0,1}\Di.\eeq
Thus, using \eqref{IndEmb2}, we see that
\eq \mc{B}=(\~n,z_{2},\dots,z_{n})\eeq
forms a basis of $T_{y}\Di$ for all $y\in \del\Di.$ 
Here we slightly abuse notation and let $\~n$ denote the contravariant version of $\~n$ as well. Furthermore we have
\eq g_{ij}\~n^{i}z^{j}_{I}=0,\quad2\leq I\leq n.\eeq
\eRem

\subsubsection*{Boundary derivatives}
\Lem\label{BoundH}
On $\del\Di$ there holds
\eq H_{i}\~n^{i}=-H. \eeq
\eLem

\pf
Note that from 
\eq \dot{X}=\fr{1}{H}N,\eeq
which also holds on $\del\Di,$ we obtain from \eqref{IndEmb2} that
\eq \fr{1}{H}x_{i}\nu^{i}=\fr{1}{H}N=\fr{d}{dt}(X\circ z)=x_{i}\dot{y}^{i},\eeq 
where $\nu$ denotes the pullback of $N$ along $x,$ which is well defined by \eqref{BoundCond}. We obtain that
\eq \label{IndFlow}\dot{y}=\fr{1}{H}\nu\eeq
holds in $T\S^{n}.$
Differentiating \eqref{BoundCond} with respect to time we obtain
\SAl 	0&=\inpr{\dot N}{\~N}+\inpr{N}{\~N_{i}\dot y^{i}}\\
		&=\fr{1}{H^{2}}\inpr{X_{i}H^{i}}{\~N}+\fr{1}{H}\inpr{N}{\check h_{i}^{k}x_{k}\nu^{i}},
\end{split}\end{align}
which implies the result in view of $\check h^{k}_{i}=\d^{k}_{i}.$
\epf

\Lem\label{Asplit}
On $\del\Di$ there hold
\begin{enumerate}[(i)]
\item{$h_{ij}\~n^{i}z^{j}_{I}=0,\quad 2\leq I\leq n,$}
\item{$h_{ij;k}z^{i}_{I}z^{j}_{J}\~n^{k}=-h_{ij}z^{i}_{I}z^{j}_{J}+h_{ij}\~n^{i}\~n^{j}g_{kl}z^{k}_{I}z^{l}_{J}.$}
\end{enumerate}
\eLem

\pf
Differentiating \eqref{BoundCond} with respect to $\xi^{I}$ yields, also using \eqref{IndEmb2},
\SAl\label{Asplit1}	0&=\inpr{N_{I}}{\~N}+\inpr{N}{\~N_{I}}\\
		&=h^{k}_{l}\~n_{k}z^{l}_{I}+\inpr{N}{\check h^{k}_{l}x_{k}y^{l}_{I}}\\
		&=h_{ij}\~n^{i}z^{j}_{I}.				\end{split}\end{align}
Differentiate \eqref{IndEmb2} twice and take the scalar product with $X_{k}$ to obtain
\eq\label{zIJ} z^{k}_{IJ}=-\check h_{lm}\~n^{k}y^{l}_{I}y^{m}_{J}=-\g_{IJ}\~n^{k}=-g_{ij}z^{i}_{I}z^{j}_{J}\~n^{k}, \eeq
where we used that $\check h_{ij}=\check g_{ij}.$

Differentiating \eqref{Asplit1} with respect to $\xi^{J}$ yields
\SAl h_{ij;k}z^{k}_{J}\~n^{i}z^{j}_{I}&=-h_{ij}\~n^{i}_{J}z^{j}_{I}-h_{ij}\~n^{i}z^{j}_{IJ}\\
						&=-h_{ij}z^{i}_{J}z^{j}_{I}+h_{ij}\~n^{i}\~n^{j}\g_{IJ}.\end{split}\end{align}
\epf

\Lem\label{ABoundary}
On $\del\Di$ there holds
\eq h_{ij;k}\~n^{i}\~n^{j}\~n^{k}=-nh_{ij}\~n^{i}\~n^{j}.\eeq
\eLem

\pf
With respect to the basis $\mc{B},$ $g$ and $A$ split, compare \cref{gsplit} and \cref{Asplit}. 
Therefore we have
\eq g^{IJ}z^{i}_{I}z^{j}_{J}=g^{ij}-\~n^{i}\~n^{j}\eeq
and thus
\SAl -H=H_{k}\~n^{k}&=g^{ij}h_{ij;k}\~n^{k}\\
		&=h_{ij;k}\~n^{i}\~n^{j}\~n^{k}+h_{ij;k}z^{i}_{I}z^{j}_{J}\~n^{k}g^{IJ}\\
		&=h_{ij;k}\~n^{i}\~n^{j}\~n^{k}-h_{ij}z^{i}_{I}z^{j}_{J}g^{IJ}+h_{ij}\~n^{i}\~n^{j}g_{kl}z^{k}_{I}z^{l}_{J}g^{IJ}\\
		&=h_{ij;k}\~n^{i}\~n^{j}\~n^{k}-H+nh_{ij}\~n^{i}\~n^{j}. \end{split}\end{align}
\epf

We need another lemma about the induced embedding.

\Lem\label{IndEmbCurv}
The second fundamental form $(\eta_{IJ})$ with respect to the normal $-\nu$ as in \eqref{IndFlow} of the induced embedding 
\eq y\cn\del \Di\hra\S^{n}\eeq
satsifies
\eq \eta_{IJ}=h_{kl}z^{k}_{I}z^{l}_{J}.\eeq
In particular, if $X$ is the embedding of a convex hypersurface into $\R^{n+1},$ $y$ is the embedding of a convex hypersurface into the sphere $\S^{n}.$
\eLem

\pf
Differentiating \eqref{IndEmb2} twice, we obtain from \eqref{zIJ}
\SAl -x_{k}\eta_{IJ}\nu^{k}&=-h_{kl}z^{k}_{I}z^{l}_{J}N+X_{k}z^{k}_{IJ}+\g_{IJ}\~N\\
					&=-h_{kl}z^{k}_{I}z^{l}_{J}N.\end{split}\end{align}
\epf

To understand how the height of our hypersurfaces over a hyperplane behaves, we have the following lemma.

\Lem\label{BoundHeight}
Let 
\eq X_{0}\cn \Di\ra M_{0}\hra\R^{n+1}\eeq
be an embedding as in \eqref{initial}. Let $\o\in\R^{n+1}.$ Then the height over the hyperplane $\o^{\perp},$
\eq w=\inpr{X}{\o},\eeq
satisfies
\eq   w_{k}\~n^{k}=w\eeq
on $\del\Di.$
In particular, if $\o$ is chosen, such that $w$ is positive on $\del\Di,$ $w$ attains its global minimum in the interior of $\Di.$
\eLem

\pf
On $\del\Di$ we have
\SAl w_{k}\~n^{k}&=\-g_{\a\b}X^{\a}_{k}\o^{\b}g^{kl}\-g_{\g\d}X^{\g}_{l}\~N^{\d}\\
			&=\-g_{\b\d}\o^{\b}\~N^{\d}\\
			&=\inpr{\~N}{\o}\\
			&=w,									\end{split}\end{align}
since on the boundary $X$ maps into $\S^{n}$ and here the position vector $X$ equals the outer normal $\~N.$
\epf

\subsubsection*{Evolution equations}
We need the following evolution equations.

\Lem\label{EvPhi}
The speed 
\eq \Phi=-\fr{1}{H}\eeq
satisfies
\eq \dot{\Phi}-\fr{1}{H^{2}}\D\Phi=\fr{\|A\|^{2}}{H^{2}}\Phi\eeq
in the interior and
\eq \Phi_{k}\~n^{k}=\Phi\eeq
on the boundary.
\eLem

\pf
The interior equation follows from \cite[Lemma~2.3.4]{Gerhardt:/2006} and the boundary derivative from \cref{BoundH}.
\epf

\Lem\label{Evw}
Let $\o\in\R^{n+1}.$ Then the height
\eq w=\inpr{X}{\o}\eeq
 of $M_{t}$ over the plane $\o^{\perp}$ satisfies
\eq \dot{w}-\fr{1}{H^{2}}\D w=\fr{2}{H}\inpr{N}{\o}\eeq
in the interior and 
\eq w_{k}\~n^{k}=w\eeq
on the boundary. 
\eLem

\pf
The interior equation comes from \eqref{Floweq} and the boundary derivative is derived in \cref{BoundHeight}.
\epf

Applying a strictly convex function in $\R^{n+1}$ to $X$ yields a very useful evolution equation, the derivation of which is a simple calculation.

\Lem\label{Evchi}
Let $\chi\in C^{2}(\R^{n+1}).$ Then $\chi=\chi(X)$ satisfies
\eq \dot{\chi}-\fr{1}{H^{2}}\D\chi=\fr{2}{H}\chi_{\a}N^{\a}-\fr{1}{H^{2}}\chi_{\a\b}X^{\a}_{i}X^{\b}_{j}g^{ij}\eeq
in the interior and 
\eq \chi_{i}\~n^{i}=\inpr{D\chi}{\~N}\eeq
on the boundary.
\eLem

\section{Height estimates}\label{Heightest}

\Def
\begin{enumerate}[(i)]
\item{For a convex hypersurface $M_{0}$ satisfying \eqref{initial} let $\mrm{conv}(\del M_{0})$ denote the convex body in the sphere enclosed by the convex hypersurface $\del M_{0}\hra\S^{n},$ cf. \cref{IndEmbCurv}.}
\item{For a point $x_{0}\in\S^{n},$ $\mc{H}(x_{0})$ denotes the closed hemisphere in $\S^{n}$ with center $x_{0}.$ The corresponding equator is denoted by $\mc{S}(x_{0}).$}
\end{enumerate}
\eDef

\Lem\label{SuppCone}
Let $M_{0}$ be a convex hypersurface satisfying \eqref{initial} and
\eq C_{0}=\{x\in\R^{n+1}\cn x=sp, s\geq 0, p\in \mrm{conv}(\del M_{0})\}.\eeq
Then there holds
\eq M_{0}\subset C_{0}.\eeq 
\eLem

\pf
$C_{0}$ is a convex cone in $\R^{n+1},$ cf. \cite[Prop.~2]{FerreiraIusemNemeth:11/2013}, and is made of an intersection of half-spaces in $\R^{n+1}$ with normal $N_{0},$
\eq C_{0}=\bigcap_{y\in \del M_{0}}\{x\in \R^{n+1}\cn \inpr{x-y}{N_{0}}\leq 0\}.\eeq
The tangent spaces of $C_{0}$ and $M_{0}$ coincide at all boundary points due to \eqref{initial3} 
and hence for all boundary points $y,$ $M_{0}$ lies on the same side of the tangent plane $T_{y}M_{0}$ as $C_{0}.$ 
\epf

In the sequel we need the following simple geometric lemma.

\Lem\label{closedCones}
Let $R>0,$ $e_{0}\in\R^{n+1}$ be a unit vector and $C\subset\R^{n+1}$ be a convex closed cone. Then for all $\e>0$ there exists $\d>0,$ such that
\eq \label{closedCones1}\inpr{a}{e_{0}}\geq \cos\(\fr{\pi}{2}-\e\)\|a\|\quad\forall a\in C\eeq
implies
\eq \inpr{x}{e_{0}}\geq R+\d\quad\forall B_{R}(x)\subset C.\eeq
\eLem

\pf
Suppose the claim was false. Then there existed $\e>0$ and a sequence of Euclidean balls $B_{R}(x_{k})\subset C$ with the property
\eq R\leq \inpr{x_{k}}{e_{0}}<R+\fr 1k\eeq
and such that \eqref{closedCones1} holds. Without loss of generality assume that $x_{k}$ converges to some $x\in C.$ Then we also have
\eq\label{closedCones2} B_{R}(x)\subset C,\eeq
since $C$ is closed. Then
\eq a=x-Re_{0}\in \-B_{R}(x)\eeq
and due to \eqref{closedCones1} there holds $a=0.$ Thus we have
\eq x=(R,0,\dots,0)\eeq
and hence a contradiction to \eqref{closedCones2}, since $C$ hits $\{x^{0}=0\}$ at $0$ transversally.
\epf

\Lem\label{Normaldirection}
Let
\eq X_{0}\cn M\hra\R^{n+1}\eeq
be the embedding of a strictly convex hypersurface $M_{0},$ such that \eqref{initial} holds. Let $e_{0}\in\mrm{int}(\mrm{conv}(\del M_{0}))$ be a direction, such that $\mrm{conv}(\del M_{0})$ is contained in the open hemisphere $\mrm{int}(\mc{H}(e_{0})).$ Then we have
\eq \p:=\inpr{N_{0}}{e_{0}}\leq C_{0}\eeq for some constant $C_{0}<0,$ which only depends on the inradius of $\mrm{conv}(\del M_{0}).$ 
\eLem

\pf
The Gauss map of the embedding $X_{0},$
\eq N_{0}\cn \Di\hra \S^{n},\eeq
is a diffeomorphism onto its image due to the strict convexity. By \cref{IndEmbCurv} and \cite[Thm.~9.2.5]{Gerhardt:/2006} the restriction 
\eq N_{0|\del\Di}\cn \S^{n-1}\hra \S^{n}\eeq
is a convex embedding and by \cite[Thm.~1.1]{CarmoWarner:/1970}, there exist two disjoint open connected components $A$ and $B,$ such that
\eq \S^{n}\backslash N_{0}(\del\Di)=A\cup B\eeq
and $A$ is the interior of the strictly convex body in the sphere, which $N_{0}(\del\Di)$ bounds.
Since $\mrm{conv}(\del M_{0})$ is chosen to be contained in $\mc{H}(e_{0}),$ we have
\eq \del A\subset \mc{H}(-e_{0})\eeq
and
from \cite[Thm.~9.2.9, Thm.~9.2.10]{Gerhardt:/2006} we obtain
\eq -e_{0}\in A\subset\-{A}\subset \mc{H}(-e_{0}).\eeq 
We have either
\eq \label{Normaldirection1}N_{0}(\Di\backslash\del\Di)\subset A\eeq
or 
\eq N_{0}(\Di\backslash\del \Di)\subset B,\eeq since the continuous map
\eq N_{0}\cn \Di\backslash\del \Di\ra A\cup B\eeq
 has to map the connected domain into a connected component, also compare \cite[Cor.~IV.19.7]{Bredon:/1993}. Since the height function
 \eq w=\inpr{X}{e_{0}}\eeq is increasing at the boundary, cf. \cref{BoundHeight}, it attains an interior minimum and thus $-e_{0}\in N_{0}(\Di\backslash\del \Di).$ Thus we must have \eqref{Normaldirection1}. This implies the claim.
\epf

\Cor\label{Maxheight}
In the situation of \cref{Normaldirection} the height function
\eq w=\inpr{X_{0}}{e_{0}}\eeq
does not attain an interior local maximum.
\eCor

\pf
Using the Gaussian formula we obtain
\eq \D w=-H\p>0.\eeq
\epf

\Cor\label{Starshaped}
In the situation of \cref{Normaldirection} there holds
\eq \inpr{e_{0}-X_{0}}{N_{0}}<0.\eeq
\eCor

\pf
Suppose the claim to be false, then there existed a point $z\in\mrm{int}(\Di)$ with the property that $e_{0}$ is not contained in the supporting open halfspace at $X_{0}=X_{0}(z),$
\eq S_{0}=\{x\in\R^{n+1}\cn\inpr{x-X_{0}}{N_{0}}< 0\}.\eeq
Due to \cref{Normaldirection} we then also had
\eq 0\notin \-{S}_{0}.\eeq
By the strict convexity of $M_{0}$ we have
\eq X_{0}(\del\Di)\subset S_{0}.\eeq
$\del S_{0}$ splits $\S^{n}$ into two spherical caps.
Translating the hyperplane $\del S_{0}$ until it hits $0,$ we see that $\del M_{0}$ originally had to be contained in the spherical cap which is geodesically convex. But by assumption we have $e_{0}\in\mrm{int}(\mrm{conv}(\del M_{0})),$ which contradicts $e_{0}\notin S_{0}.$  
\epf

We are now able to estimate the height of a hypersurface $M_{0}$ as the latter appears in \eqref{initial}. It depends on the estimate in \cref{Normaldirection} and the curvature.

\Lem\label{PositiveHeight}
In the situation of \cref{Normaldirection} the height
\eq w=\inpr{X}{e_{0}}\eeq
satisfies
\eq w\geq \d>0,\eeq
for a constant $\d,$ which depends on the constant $C_{0}$ in \cref{Normaldirection}, the length of the second fundamental form of $M_{0}$ and the distance of $\del M_{0}$ to the equator $\mc{S}(e_{0}).$
\eLem

\pf
Let $a\in M_{0}$ be the interior global minimum point of $w.$ Due to \cref{Normaldirection} it is possible to write $M_{0}$ locally around $a$ as a graph over the unit disk in $\{0\}\times\R^{n},$ where $w$ is the graph function. Then
\eq w_{ij}=-h_{ij}\inpr{N}{e_{0}}.\eeq
Using \cite[Lemma~2.7.6]{Gerhardt:/2006}, we obtain that the Hessian of $w$ with respect to Euclidean coordinates only depends on the second fundamental form and on the estimate of $\inpr{N}{e_{0}}$ from below. Define
\eq \hat{M}_{0}=\bigcap_{y\in M_{0}}\{x\in \R^{n+1}\cn \inpr{x-y}{N_{0}}\leq 0\}.\eeq
From the previous considerations $\hat M_{0}$ satisfies an interior sphere condition at $a$ with interior ball $B_{R}$ depending on $\sup\|A\|$ and $\inpr{N}{e_{0}}.$  
 Due to 
\eq B_{R}\subset\hat{M}_{0}\subset C_{0},\eeq 
from \cref{closedCones} we obtain the existence of $\d>0,$ such that 
\eq \inpr{a}{e_{0}}\geq \d.\eeq
\epf

\Cor\label{InsideUnitBall}
In the situation of \cref{Normaldirection} we have \eq X_{0}(\mrm{int}(\Di))\subset \mrm{int}(B^{+}),\eeq
where $B^{+}\subset\R^{n+1}$ is the pointed halfball
\eq B^{+}=B_{1}^{+}(0)\backslash\{e_{0}\}.\eeq
\eCor

\pf
The function
\eq \rho=|X_{0}|^{2}\eeq
satisfies
\eq \D\rho=-2H\inpr{N_{0}}{X_{0}}+2n,\eeq
due to the Gaussian formula.
At an interior maximum of $\rho$ we have
\eq 0=\nabla\rho\eeq
and thus $X_{0}$ has to be a multiple of $N_{0}.$ Since
\eq \inpr{X_{0}}{e_{0}}>0\eeq
due to \cref{PositiveHeight} and 
\eq \inpr{N_{0}}{e_{0}}<0\eeq
due to \cref{Normaldirection}, we have
\eq \inpr{N_{0}}{X_{0}}<0.\eeq
Thus at a maximal point we have
\eq \D\rho >0,\eeq
a contradiction. Since we have $\rho=1$ at the boundary, the claim follows. 
\epf

\section{Moebius coordinates and the scalar flow}\label{MoebiusCoord}
In this section we want to derive a scalar flow equation naturally associated with \eqref{Flow}. Therefore we aim for a graph representation. A natural candidate for hypersurfaces of our type are rotations of Moebius transformations on the plane. Consider a one-parameter familiy of Moebius transformations of the form
\eq \~f(x,\l)=\fr{(1+\l)x+i(\l-1)}{1+\l+i(1-\l)x},\eeq
where $(x,\l)\in [-1,1]\times [1,\infty).$ For each $\l$ this is a conformal transformation moving the real axis towards $i$ as $\l\ra\infty,$ whereas the boundary of the real interval $[-1,1]$ maps to the unit sphere perpendicularly. A rotation of a plane in $\R^{n+1}$ around the $e_{0}$-axis gives rise to the following definition.

\Def
Let $D\subset \R^{n}$ be the unit disk. Define \textit{Moebius coordinates} for the pointed halfball 
\eq B^{+}:=B_{1}^{+}(0)\backslash\{e_{0}\}\eeq
to be the diffeomorphism
\SAl\label{MoebiusDiffeo}	f&\cn D\times [1,\infty)\ra B^{+}\\
		f(x,\l)&=\fr{4\l x+(1+|x|^{2})(\l^{2}-1)e_{0}}{(1+\l)^{2}+(1-\l)^{2}|x|^{2}}.\end{split}\end{align}
\eDef 

\subsubsection*{Graphs in Moebius coordinates}
Let us provide some general formalae for hypersurfaces $M\subset\R^{n+1}$ which can be written as graphs in Moebius coordinates. Thus suppose the embedding of a hypersurface $M$ is given by a map
\SAl X\cn \Di&\hra\R^{n+1}\\
		z&\mt f(x(z),u(x(z))),\end{split}\end{align}
where $u\cn D\ra[1,\infty)$ is a function. First of all, from a tedious computation and the conformality of $f$ we obtain a representation of the Euclidean metric $\d_{\a\b}$ in Moebius coordinates,
\eq d\-{s}^{2}=e^{2\psi}({dx^{0}}^{2}+\s_{ij}dx^{i}dx^{j}),\eeq
where $x^{0}$ corresponds to the $\l$-coordinate,
\eq e^{2\psi}=\left\langle \fr{\del f}{\del x^{0}},\fr{\del f}{\del x^{0}}\right\rangle,\eeq
\eq\label{TangVecf} \fr{\del f}{\del \l}(x,\l)=\fr{(1+|x|^{2})(1-\l^{2})}{\l((1+\l)^{2}+(1-\l)^{2}|x|^{2})}\(f-\fr{\l^{2}+1}{\l^{2}-1}e_{0}\).\eeq
and
\eq \s_{ij}=e^{-2\psi}\left\langle \fr{\del f}{\del x^{i}},\fr{\del f}{\del x^{j}}\right\rangle.\eeq
For $M$ we have the induced metric
\eq g_{ij}=e^{2\psi}(u_{i}u_{j}+\s_{ij})\eeq
with inverse
\eq g^{ij}=e^{-2\psi}\(\s^{ij}-\fr{\s^{ik}u_{k}}{v}\fr{\s^{lj}u_{l}}{v}\),\eeq
where
\eq v^{2}=1+\s^{ij}u_{i}u_{j}.\eeq
The contravariant version of the normal is
\eq (N^{\a})=\pm v^{-1}e^{-\psi}(1,-\s^{ik}u_{k}).\eeq
Those formulae can be found in \cite[Sec.~1.5]{Gerhardt:/2006}.

Due to the conformality of $f$ the outward Euclidean unit normal to $D,$ $\breve N,$ is mapped to a multiple of the unit normal to the sphere in $\R^{n+1}$ which we called $\~N$ earlier. Thus for a hypersurface satisfying the boundary condition \eqref{initial3} we obtain
\SAl	0&=\left\langle \breve N^{k}\fr{\del f}{\del x^{k}},N\right\rangle\\
		&=\mp \fr{e^{\psi}}{v}\breve N^{k}u_{k}\end{split}\end{align}  
and thus such a hypersurface satisfies the Neumann boundary condition
\eq \label{NeumannCondition}\breve N^{k}u_{k}=0.\eeq

Now we prove that hypersurfaces satisfying \eqref{initial} are graphs in Moebius coordinates. 

\Prop\label{Moebius}
Let 
\eq X_{0}\cn M\hra\R^{n+1}\eeq be the embedding of a strictly convex hypersurface $M_{0},$ such that \eqref{initial} holds. Choose $e_{0}\in\mrm{int}(\mrm{conv}(\del M_{0})),$ such that $\mrm{conv}(\del M_{0})$ is contained in the open hemisphere $\mc{H}(e_{0}).$ Then $M_{0}$ can be written as a graph in Moebius coordinates around $e_{0},$ i.e. Moebius coordinates in the pointed half-ball $B_{1}^{+}(0)\backslash\{e_{0}\}$ yield a representation 
\eq X_{0}(z)=f(x,u_{0}(x)),\eeq
where $f$ is the diffeomorphism defined in \eqref{MoebiusDiffeo}.
\eProp

\pf
Due to \cref{InsideUnitBall} Moebius coordinates are well-defined throughout $M_{0}.$
By the implicit function theorem all we have to show is that
\eq \left\langle \fr{\del f}{\del \l},N_{0}\right\rangle<0.\eeq
Due to \cref{PositiveHeight} we have $\l\geq c>1$ and thus it suffices to discard the negative scalar fraction in \eqref{TangVecf}. We have
\SAl \label{Moebius1}\left\langle X_{0}-\fr{\l^{2}+1}{\l^{2}-1}e_{0},N_{0}\right\rangle&=\left\langle X_{0}-e_{0},N_{0}\right\rangle-\left\langle \fr{2}{\l^{2}-1}e_{0},N_{0}\right\rangle\\	
											&>-\fr{2}{\l^{2}-1}\inpr{e_{0}}{N_{0}}			\\
											&>0,			\end{split}\end{align}
due to \cref{Normaldirection} and \cref{Starshaped}.
\epf

The previous considerations allow us to naturally associate a scalar parabolic equation to strictly convex solutions of our inverse mean curvature flow \eqref{Flow}.

\Cor\label{ScalarProblem}
Let $X$ be a solution of \eqref{Flow} on a time interval $[0,\e),$ such that all $M_{t},$ $0\leq t<\e,$ range within a pointed halfball $B^{+}$ and are graphs in Moebius coordinates for $B^{+},$
\eq M_{t}=\{(x(t,z),u(t,x))\cn (t,z)\in [0,\e)\times \Di\}.\eeq
Then $u$ solves a parabolic Neumann problem on $[0,\e)\times D,$ namely
\PBVP{\fr{\del u}{\del t}}{-\fr{v}{e^{\psi}H}}{(0,\e)\times D,}{u_{k}\breve N^{k}}{0}{[0,\e)\times \del D,}{u}{u_{0}}{\{0\}\times D.}
\eCor

\pf
For curvature flows in ambient spaces covered by Gaussian coordinate systems the interior equations are deduced in \cite[p.~98-99]{Gerhardt:/2006}. Just note that in our case the normal $N_{0}$ and the vector $\fr{\del f}{\del x^{0}}$ are pointing in opposite directions, hence the sign.
The boundary equation follows from the fact that all $M_{t}$ are perpendicular to the sphere and by the derivation of \eqref{NeumannCondition}.
\epf

\section{Curvature estimates and convexity}\label{CurvEst}
\Rem
Let $T^{*}$ be the largest time, such that there exists a smooth solution to \eqref{Flow} on the interval $[0,T^{*}).$ This implies mean convexity of $M_{t},$ $0\leq t< T^{*}.$
By \cref{Shorttime} we indeed have $T^{*}>0.$
Let $\-T>0$ be the largest time, such that the solution is smooth on $[0,\-T)$ and $M_{t}$ is strictly convex for all $0\leq t<\-T.$
\eRem

\Prop\label{kappaBound}
Let $X$ be the solution of \eqref{Flow} on the interval $[0,\-T).$ Then the principal curvatures are bounded, i.e. for $1\leq i\leq n$ there holds 
\eq \k_{i}\leq H\leq \max_{\Di}H(0,\cdot)\quad \forall t\in[0,\-T).\eeq
\eProp

\pf
Using the convexity of the flow hypersurfaces up to $\-T,$ all we have to bound is $H.$ From \cref{EvPhi} we obtain
\eq \dot{H}-\frac{1}{H^{2}}\D H\leq -\fr{\|A\|^{2}}{H^{2}}H\eeq
 and \eq H_{k}\~n^{k}=-H.\eeq
 Thus the claim follows from a standard maximum principle, e.g. \cite[Thm.~3.1]{Stahl:06/1996}.
\epf

\Lem
On the interval $[0,\-T)$ let 
\eq y_{t}\cn \del\Di\hra\S^{n}\eeq
be the induced embeddings of $X_{t}.$ Then the convex bodies of the embedded submanifolds $\del M_{t}\hra \S^{n}$ form an increasing sequence and satisfy uniform interior sphere conditions independently of $t.$
\eLem

\pf
The convexity of the $\del M_{t}$ in $\S^{n}$ follow from \cref{IndEmbCurv}. From \eqref{IndFlow} we see that the enclosed convex bodies are increasing. From \cref{kappaBound} and \cref{IndEmbCurv} we obtain uniform $C^{2}$-estimates and thus uniform interior sphere conditions, also compare \cite[Def.~3.2]{MakowskiScheuer:/2013}.
\epf

\Cor
There exists a $C^{1,\a}$ limiting surface $\del M_{\-T}$ arising as the limit of the $\del M_{t}.$ $\del M_{\-T}$ either is an equator of the sphere or is containd in an open hemisphere.
\eCor

\pf
$\del M_{\-T}$ is the boundary of a weakly convex body in a hemisphere, in the sense of \cite[Def.~3.2]{MakowskiScheuer:/2013}, also compare \cite[Lemma~6.1]{MakowskiScheuer:/2013}. \cite[Thm.~1.1]{MakowskiScheuer:/2013} implies the claim.
\epf

We want to conclude that $\-T=T^{*}$ and that $\del M_{T^{*}}$ must be an equator, which would yield the result due to the height estimates. Therefore we need some more estimates.

\Lem\label{Hbelow}
Let $X$ be the solution of \eqref{Flow} on the interval $[0,\-T)$ and suppose that $\del M_{\-T}$ is not an equator. Then there holds
\eq \sup\limits_{[0,\-T)\times \Di}\fr{1}{H}\leq c,\eeq
where $c$ depends on $M_{0}$ and the distance of $\del M_{\-T}$ to a suitable equator $\mc{S}(e_{0}).$
\eLem

\pf
Let $e_{0}\in\mrm{int}(\mrm{conv}(\del M_{\-T})),$ such that $\mrm{conv}(\del M_{\-T})$ is contained in $\mrm{int}(\mc{H}(e_{0})).$ Then, due to the monotonicity of $\mrm{conv}(\del M_{t})$ we also have
\eq e_{0}\in\mrm{int}(\mrm{conv}(\del M_{t}))\eeq
for $t$ close to $\-T.$
Thus it is possible to apply \cref{PositiveHeight} to obtain a positive lower bound for the height function
\eq w=\inpr{X_{t}}{e_{0}}\geq \d>0.\eeq
Define the strictly convex function in $\R^{n+1}$
\eq \chi(x)=\fr{1}{2}|\hat x|^{2}+\fr{\b}{2}(x^{0})^{2}-\l x^{0}+1,\eeq
where 
\eq \hat x=(0,x^{1},\dots,x^{n})\eeq
and
\eq \l>\fr{1}{\d},\quad 0<\b<1.\eeq
Define
\eq \z=\fr{1}{H}\fr{1}{\fr{1}{2}-\chi}\equiv\fr{1}{H}G(\chi).\eeq
Due to the height estimates, $\z$ is well defined and positive on $[0,\-T)\times \Di.$
With the help of \cref{EvPhi} and \cref{Evchi} a simple computation yields the following evolution equation for $\z,$ namely
\SAl\label{Hbelow1}	\dot{\z}-\fr{1}{H^{2}}\D\z&=\fr{\|A\|^{2}}{H^{2}}\z+2\chi_{\a}N^{\a}\z^{2}-\fr{1}{H}\chi_{\a\b}X^{\a}_{i}X^{\b}_{j}g^{ij}\z^{2}\\
					&\hp{=}-2\chi_{i}\chi^{i}\z^{3}-\fr{2}{H^{2}}\(\fr{1}{H}\)_{i}G^{i}\end{split}\end{align}
and the boundary equation
\eq \z_{i}\~n^{i}=\(1+G\chi_{\a}\~N^{\a}\)\z.\eeq
Due to $X=\~N$ on the boundary, we obtain
\eq \chi_{\a}\~N^{\a}=1+(\b-1)(X^{0})^{2}-\l X^{0}\eeq
and thus on the boundary
\eq 1+G\chi_{\a}\~N^{\a}=1+\fr{1+(\b-1)(X^{0})^{2}-\l X^{0}}{\l X^{0}-\fr{\b-1}{2}(X^{0})^{2}-1}<0.\eeq

Now suppose for $0<T<\-T$ that
\eq \max\limits_{[0,T]\times\Di}\z=\z(t_{0},z_{0})\geq 1,\quad t_{0}>0.\eeq
Then $z_{0}\in \mrm{int}(\Di)$ and thus from \eqref{Hbelow1} we obtain at this point that, also using
\eq \fr{G_{i}}{G}=-\fr{\(\fr{1}{H}\)_{i}}{\fr 1H},\eeq
\eq 0\leq \(c-\fr{1}{H}\chi_{\a\b}X^{\a}_{i}X^{\b}_{j}g^{ij}\)\z^{2},\eeq
where $c=c(\d).$ Since
\SAl \chi_{\a\b}X^{\a}_{i}X^{\b}_{j}g^{ij}&=\chi_{\a\b}\-g^{\a\b}-\chi_{\a\b}N^{\a}N^{\b}\\
							&=n+\b-1+(1-\b)(N^{0})^{2},				\end{split}\end{align}
we obtain a bound for $\fr{1}{H}$ at the point $(t_{0},z_{0}).$ Since $G$ is bounded, this implies a uniform bound on $\z$ and in turn a uniform bound on $\fr{1}{H}.$
\epf

\Prop\label{PresConv}
There holds $\-T=T^{*}.$ In particular the strict convexity of the flow hypersurfaces is preserved up to $T^{*}.$
\eProp

\pf
Suppose that $\-T<T^{*}\leq\infty.$ In case that $\del M_{\-T}$ is an equator of the sphere, we conclude from the height estimates that $M_{\-T}$ is a flat disk and thus a singularity of the flow. This would yield $\-T=T^{*}.$ Thus suppose that $\del M_{\-T}$ is not an equator. From \cref{Hbelow} we obtain 
\eq \fr 1H\leq c\quad\forall t\in[0,\-T)\eeq
and again the height function satisfies
\eq w\geq \d>0.\eeq
Define
\eq \~H=\sum_{i=1}^{n}\fr{1}{\k_{i}}=g_{ij}\~h^{ij},\eeq
where $(\~h^{ij})$ is the inverse of $(h_{ij}).$
At a given point choose coordinates with respect to the basis $\mc{B}=(\~n,z_{I}),$ then at the boundary we deduce, due to \cref{Asplit} and \cref{ABoundary}, that
\SAl \~H_{k}\~n^{k}&=-\~h^{r}_{i}\~h^{si}h_{rs;k}\~n^{k}\\
				&=-\~h^{1}_{1}\~h^{11}h_{11;k}\~n^{k}-\~h^{J}_{I}\~h^{KI}h_{JK;k}\~n^{k}\\
				&=n\~h^{1}_{1}\~h^{11}h_{11}+\~h^{J}_{I}\~h^{KI}h_{JK}-\~h^{J}_{I}\~h^{KI}g_{KJ}h_{11}\\
				&\leq (n-1)\~h^{1}_{1}+\~h^{r}_{i}\~h^{si}h_{rs}\\
				&=(n-1)\~h^{i}_{j}\~n_{i}\~n^{j}+\~H.		\end{split}\end{align}

Set \eq \phi=\log \~H-(n+1)\log w-\a t,\quad t<\-T,\eeq
where 
$\a$ will be chosen in dependence of $\d$ and the initial data.
From \cite[Lemma~6.5]{Gerhardt:01/1996} and \cref{Evw} we obtain
\SAl	\label{PresConv1}	\dot{\phi}-\fr{1}{H^{2}}\D\phi&=-\fr{\|A\|^{2}}{H^{2}}+\fr{2n}{H\~H}+\fr{2}{H^{2}\~H^{2}}\~H_{i}\~H^{i}\\
							&\hp{=}-\(\fr{2}{H^{2}}g^{rs}\~h^{kl}h_{rk;p}h_{sl;q}-\fr{2}{H^{3}}H_{p}H_{q}\)\fr{\~h^{pi}\~h^{q}_{i}}{\~H}\\
							&\hp{=}-\fr{2n+2}{Hw}\inpr{N}{e_{0}}-\fr{n+1}{H^{2}w^{2}}w^{i}w_{i}-\a		 \end{split}\end{align}
in the interior and
\eq \phi_{k}\~n^{k}\leq 1+\fr{n-1}{\~H}\~h^{1}_{1}-(n+1)<-1\quad\forall (t,\xi)\in[0,\-T)\times \del \Di.\eeq
Now suppose that for $0<T<\-T$ we have
\eq \sup\limits_{[0,T]\times \Di}\phi=\phi(t_{0},z_{0}),\quad t_{0}>0.\eeq
Then $z_{0}$ does not lie on  $\del\Di.$ From \eqref{PresConv1} we obtain at $(t_{0},z_{0}),$
also using
\eq \fr{\~H_{i}}{\~H}=(n+1)\fr{w_{i}}{w}\eeq
and that the big bracket is nonnegative by \cite[equ.~(1.7)]{Gerhardt:01/1996}, 
 that
\eq 0\leq c+\fr{2(n+1)^{2}}{H^{2}w^{2}}\|Dw\|^{2}-\a,\eeq
where the constant depends on $\d$ and the bound on $H^{-1}.$
For large $\a$ this is a contradiction. Thus under the assumption that $\del M_{\-T}$ is not an equator we obtain that the supremum of $\phi$ would be decreasing and thus $\phi$ was bounded up to $\-T.$ But then
\eq \log\~H=\phi+(n+1)\log w+\a t\leq c+\a\-T,\eeq
which contradicts the definition of $\-T,$ at which $\~H$ would have to blow up, provided $\-T<T^{*}.$ 
\epf

\Cor\label{MaxExist}
There holds
\eq T^{*}<\infty.\eeq
\eCor

\pf
Let $e_{0}\in\mrm{int}(\mrm{conv}(\del M_{T^{*}})),$ such that $\mrm{conv}(\del M_{T^{*}})\subset {\mc{H}(e_{0})}.$ The induced strictly convex hypersurfaces $\del M_{t}\hra\S^{n}$ satisfy the flow equation \eqref{IndFlow}, which has a uniformly positive speed in normal direction. Thus $\del M_{T^{*}}$ is reached in finite time.
\epf

\section{Convergence to a flat disk}\label{Flatdisk}
We have seen that as long as the boundary of the flow is strictly contained in an open hemisphere, we have uniform bounds on the height, the mean curvature and the principal curvatures. We want to conclude that the flow can be extended whenever $\del M_{T^{*}}$ is not an equator. This would finish the proof of the main result due to the definition of $T^{*}.$
In this section we will apply regularity theory to the scalar flow equation in \cref{ScalarProblem} to achieve this.

A straightforward computation yields the following representation of this equation.

\Prop\label{UniParabolic}
The function $u\cn (0,T^{*})\times D\ra [1,\infty)$ satisfies the equation
\eq \fr{\del u}{\del t}=-\fr{v}{e^{2\psi} v^{-1}g^{ij}u_{i,j}+A(x,u,Du)}\equiv F(x,u,Du,D^{2}u),\eeq
where $A$ is smooth and $F$ is a uniformly parabolic operator, provided $\del M_{T^{*}}$ is not an equator of the sphere. 
\eProp

\pf
An easy computation gives a relation between covariant and partial derivatives of $u,$ namely
\eq u_{ij}=u_{i,j}v^{-2}+r_{ij}(x,u,Du),\eeq
where $r_{ij}$ is a smooth tensor of the indicated variables.
Due to \cite[equ.~(1.5.10)]{Gerhardt:/2006} we obtain
\eq \label{UniParabolic1}h_{ij}v^{-1}\psi^{-1}=u_{i,j}v^{-2}+r_{ij}(x,u,Du)\eeq
with a possibly different, but still smooth, tensor $r_{ij}.$ Inserting this into \eqref{PBVP} gives the first equality.

The parabolicity follows from
\eq \fr{\del F}{\del u_{i,j}}=\fr{v}{e^{\psi} H^{2}}\fr{\del H}{\del u_{i,j}}=\fr{1}{H^{2}}g^{ij},\eeq
since as long as $\del M_{T^{*}}$ is not an equator, we have $H\geq c>0$ by \cref{Hbelow} and $g^{ij}$ is equivalent to the Euclidean metric on $D$ due to \eqref{Moebius1}.
\epf

\Lem
Let $X\cn(0,T]\ra \R^{n+1}$ be a solution of \eqref{Flow} and suppose that $\del M_{T}$ is not an equator of the sphere. Then
\eq T^{*}>T+\e,\eeq
where $\e$ depends on $M_{0}$ and the distance of $\del M_{T}$ to a suitable equator.
\eLem

\pf
(i) Considering the scalar proplem as in \cref{ScalarProblem}, from \cref{UniParabolic} and standard regularity theory we obtain $C^{\infty}$-estimates up to $T$ for $u,$ compare for example \cite[Thm.~14.23]{Lieberman:/1998} or \cite[Thm.~4, Thm.~5]{Uraltseva:06/1994}. A slight adjustment of the proof of \cite[Thm.~2.5.7]{Gerhardt:/2006} to the Neumann boundary case\footnote{Also compare \href{http://www.math.uni-heidelberg.de/studinfo/gerhardt/KleineZeitenNeumann.pdf}{www.math.uni-heidelberg.de/studinfo/gerhardt/KleineZeitenNeumann.pdf}} yields a short-time existence interval of length $\e$ for $C^{2,\a}$ initial functions, depending on the data of the differential operator. In our situation, these data are uniformly under control, such that choosing a flow hypersurface $M_{t_{0}}$ with $T-t_{0}<\e$ yields an extension beyond $T.$ By the standard method of difference quotients this extension is smooth. Thus we have extended the scalar function $u.$

(ii) To obtain the full curvature flow from the scalar function $u,$ we use the standard method applied in \cite[Sec.~2.3]{Marquardt:/2012}, solving an ODE to allow for normal directed evolution.
\epf

Together with \cref{MaxExist} and the $C^{2}$-estimates we obtain the final result.

\Cor
$\del M_{T^{*}}$ is an equator of the sphere and $M_{T^{*}}$ is an embedded flat disk.
\eCor

\Rem\label{Convergence}
From \cref{kappaBound} and \eqref{UniParabolic1} we obtain uniform $C^{2}$-bounds for the graph functions $u$ and thus the norm of convergence, in which the flow hypersurfaces converge to unit disk can be characterized by saying that the functions $u$ converge to the constant function with value $1$ in the norm of $C^{1,\b}(D).$
\eRem

\subsection*{Acknowledgements}
We would like to thank Oliver Schn{\"u}rer for suggesting the investigation of this interesting problem.

\bibliographystyle{/Users/J_Mac/Documents/Uni/TexTemplates/hamsplain}
\bibliography{/Users/J_Mac/Documents/Uni/TexTemplates/Bibliography}

\providecommand{\bysame}{\leavevmode\hbox to3em{\hrulefill}\thinspace}
\providecommand{\href}[2]{#2}
\begin{thebibliography}{10}

\bibitem{Bredon:/1993}
Glen Bredon, \emph{Topology and {G}eometry}, Graduate Texts in Mathematics,
  Springer, 1993.

\bibitem{CarmoWarner:/1970}
Manfredo do~Carmo and Frank Warner, \emph{Rigidity and convexity of
  hypersurfaces in spheres}, J. Differ. Geom. \textbf{4} (1970), 133--144.

\bibitem{FerreiraIusemNemeth:11/2013}
Orizon Ferreira, Alfredo Iusem, and Sandor Nemeth, \emph{Projections onto
  convex sets on the sphere}, J. Glob. Optim. \textbf{57} (2013), no.~3,
  663--676.

\bibitem{Gerhardt:/1990}
Claus Gerhardt, \emph{Flow of nonconvex hypersurfaces into spheres}, J. Differ.
  Geom. \textbf{32} (1990), no.~1, 299--314.

\bibitem{Gerhardt:01/1996}
\bysame, \emph{Closed {W}eingarten hypersurfaces in space forms}, Geom. Anal.
  Calc. Var. (1996), 71--98.

\bibitem{Gerhardt:/2006}
\bysame, \emph{Curvature problems}, Series in Geometry and Topology,
  vol.~39, International Press of Boston Inc., 2006.

\bibitem{Gerhardt:11/2011}
\bysame, \emph{Inverse curvature flows in hyperbolic space}, J. Differ.
  Geom. \textbf{89} (2011), no.~3, 487--527.

\bibitem{GigaOhnumaSato:05/1999}
Yoshikazu Giga, Masaki Ohnuma, and Moto-Hiko Sato, \emph{On the strong maximum
  principle and the large time behavior of generalized mean curvature flow with
  the {N}eumann boundary condition}, J. Differ. Equ. \textbf{154} (1999),
  no.~1, 107--131.

\bibitem{Huisken:02/1989}
Gerhard Huisken, \emph{Non-parametric mean curvature evolution with boundary
  conditions}, J. Differ. Equ. \textbf{77} (1989), no.~2, 369--378.

\bibitem{HuiskenIlmanen:/2001}
Gerhard Huisken and Tom Ilmanen, \emph{The inverse mean curvature flow and the
  {R}iemannian {P}enrose inequality}, J. Differ. Geom. \textbf{59} (2001),
  no.~3, 353--437.

\bibitem{Lambert:/2014}
Ben Lambert, \emph{The constant angle problem for mean curvature flow inside
  rotational tori}, to appear in Math. Res. Lett. \textbf{21}, no. 3, 2014,
  {\href{http://arxiv.org/abs/1207.4422}{arxiv:1207.4422}}.

\bibitem{Lieberman:/1998}
Gary Lieberman, \emph{Second order parabolic differential equations}, World
  Scientific, 1998.

\bibitem{MakowskiScheuer:/2013}
Matthias Makowski and Julian Scheuer, \emph{Rigidity results, inverse curvature
  flows and {A}lexandrov-{F}enchel type inequalities in the sphere}, preprint,
  2013, \href{http://arxiv.org/abs/1307.5764}{arxiv:1307.5764v2}.

\bibitem{Marquardt:/2012}
Thomas Marquardt, \emph{The inverse mean curvature flow for hypersurfaces with
  boundary}, Ph.D. thesis, Freie Universit{\"a}t Berlin, 2012.

\bibitem{Marquardt:07/2013}
\bysame, \emph{Inverse mean curvature flow for star-shaped hypersurfaces
  evolving in a cone}, J. Geom. Anal. \textbf{23} (2013), no.~3, 1303--1313.

\bibitem{Scheuer:06/2014}
Julian Scheuer, \emph{Non-scale-invariant inverse curvature flows in hyperbolic
  space}, Calc. Var. Partial Differ. Equ. (2014), online,
  {\href{http://link.springer.com/article/10.1007%2Fs00526-014-0742-9}{doi:10.1007/s00526-014-0742-9}}.

\bibitem{Stahl:06/1996}
Axel Stahl, \emph{Regularity estimates for solutions to the mean curvature flow
  with a {N}eumann boundary condition}, Calc. Var. Partial Differ. Equ.
  \textbf{4} (1996), no.~4, 385--407.

\bibitem{Stahl:08/1996}
\bysame, \emph{Convergence of solution to the mean curvature flow with a
  {N}eumann boundary condition}, Calc. Var. Partial Differ. Equ. \textbf{4}
  (1996), no.~5, 421--441.



\bibitem{Uraltseva:06/1994}
Nina Ural'tseva, \emph{A nonlinear problem with an oblique derivative for
  parabolic equations}, J. Math. Sci. \textbf{70} (1994), no.~3, 1817--1827.

\bibitem{Urbas:/1990}
John Urbas, \emph{On the expansion of starshaped hypersurfaces by symmetric
  functions of their principal curvatures}, Math. Z. \textbf{205} (1990),
  no.~1, 355--372.

\end{thebibliography}

\end{document}